\documentclass[12pt]{amsart}

\usepackage{epsfig}

\renewcommand{\subjclassname}{\textup{2000} Mathematics Subject
     Classification}

\newtheorem{theorem}{Theorem}[section]
\newtheorem{lemma}[theorem]{Lemma}
\newtheorem{proposition}[theorem]{Proposition}
\newtheorem{corollary}[theorem]{Corollary}

\theoremstyle{definition}

\newtheorem{example}[theorem]{Example}

\theoremstyle{remark}
\newtheorem{remark}[theorem]{Remark}

\numberwithin{equation}{section}

\def\R{{\mathbb R}}

\def\N{{\mathbb N}}
\def\Z{{\mathbb Z}}
\def\D{{\mathbb D}}
\def\Q{{\mathbb Q}}

\def\C{{\mathbb C}}

\def\cF{{\mathcal F}}

\def\cS{{\mathcal S}}

\def\Re{{\rm Re}\,}
\def\Im{{\rm Im}\,}
\def\TT{(T(t))_{t\geq 0}}

\def\TcT{({\mathbf T} (t))_{t\geq 0}}

\def\TcT0{({\mathbf T}_0 (t))_{t\geq 0}}

\def\L1{L^1 (\R_+ )}

\begin{document}

\title[Uniqueness theorems with applications]{Uniqueness theorems for (sub-)harmonic functions with applications
to operator theory}

\author{Alexander Borichev}
\address{LABAG, Universit\'e Bordeaux 1, 351, cours de la Lib\'eration, 33405 Talence cedex, France}
\email{Alexander.Borichev@math.u-bordeaux1.fr}

\author{Ralph Chill}
\address{Laboratoire de Math\'ematiques et Applications de Metz, UMR 7122, Universit\'e de Metz et CNRS, B\^at. A, Ile du Saulcy, 57045 Metz Cedex 1, France}
\email{chill@univ-metz.fr}

\author{Yuri Tomilov}
\address{Faculty of Mathematics and Computer Science, Nicolas Copernicus University, ul. Chopina 12/18, 87-100 Torun, Poland}
\email{tomilov@mat.uni.torun.pl}
\thanks{The third author was partially supported by a KBN grant.}


\date{December 10, 2005}



\begin{abstract}
We obtain uniqueness theorems for harmonic and subharmonic
functions of a new type. They lead to new analytic extension
criteria and new conditions for  stability of operator semigroups
in Banach spaces with Fourier type.
\end{abstract}

\renewcommand{\subjclassname}{\textup{2000} Mathematics Subject Classification}

\maketitle

\section{Introduction}


If $u$ is a non-negative subharmonic function on the unit disc
$\D$ such that $\lim_{z \to \xi}u(z) = 0$ for every $\xi$ on the
unit circle ${\mathbb T}$, then the classical maximum principle forces $u$ to be zero. In this article we shall present several generalisations of this uniqueness principle. 

It is natural to ask whether the above uniqueness result remains true under weaker conditions on the boundary behaviour of $u$. For example, we may consider only nontangential or even radial boundary values. However, in this case the uniqueness principle is no longer true unless 
we impose some additional restrictions on $u$; if $P(r,\varphi)$ is the Poisson kernel in $\D$, 
then the function
\begin{equation}\label{u_0}
u(re^{i \varphi})=\frac{\partial P(r, \varphi)}{\partial \varphi}= -
\sum_{n=1}^{\infty} n r^n \sin n\varphi,
\end{equation}
is harmonic in $\D$, $\lim_{r \to 1-}
u(re^{i\varphi})=0$ for every $\varphi \in [0, 2\pi]$, but clearly we have $u\ne 0$, \cite{Sh63}, \cite{Da77}, \cite{AsBr91}. In the
class of subharmonic functions for which
$$
M_r (u):=\sup_{\varphi\in [0,2\pi )} | u(re^{i\varphi } ) |, \quad
r \in [0,1),
$$
grows sufficiently slowly as $r \to 1-$
it does suffice to consider only radial boundary values;
the following radial uniqueness theorem was proved by B. Dahl\-berg, \cite[Theorem 1]{Da77}.

\begin{theorem}[Dahlberg]\label{dal}
Let $u$ be a non-negative subharmonic function on the unit disc. Assume that
\begin{eqnarray*}
(i) & & \lim_{r \to 1-} u(re^{i\varphi}) = 0 \mbox{ for every } \varphi\in [0,2\pi );\\
(ii) & & M_r (u) = {\rm o} ( (1-r)^{-2} ) \quad \text{ as } r\to 1- .
\end{eqnarray*}
Then $u=0$.
\end{theorem}

The second condition cannot be improved by replacing ``o'' by ``O'' as the example from \eqref{u_0} shows. Thus, Theorem \ref{dal} is optimal with respect to the growth conditions on $M_r(u)$.

If $u(z)$ tends to zero as $z \to \xi$ non-tangentially (that is, $z$ stays in a sector with 
a fixed opening and vertex at $\xi$)
for every $\xi \in {\mathbb T}$, then the growth condition on $M_r (u)$ in Theorem \ref{dal} can be relaxed, depending on the opening of the
sector. Such sectorial variant of Theorem \ref{dal} was obtained
by R. Berman and W. Cohn, \cite[Theorem 1]{BeCo88}. For $\theta\in (0,\frac{\pi}{2})$ and $\varphi\in [0,2\pi )$, define the sector $\Omega_{\theta}(\varphi):=\{z \in \D : |\arg (1- e^{-i\varphi}z) | \le \frac{\pi}{2} - \theta \}$.

\begin{theorem}[Berman-Cohn] \label{beco}
Let $u$ be a non-negative subharmonic function on the unit disc, and
let $\theta \in (0,\frac{\pi}{2})$ be fixed. Assume
that
\begin{eqnarray*}
(i) & & \lim_{z\to e^{i\varphi} \atop z \in \Omega_{\theta}(\varphi)}
 u(z) = 0 \mbox{ for every } \varphi \in [0,2\pi );\\
(ii) & & M_r (u) = {\rm o} ( (1-r)^{-\frac{\pi}{\theta}} ) \quad \text{ as } r\to 1- .
\end{eqnarray*}
Then $u=0$.
\end{theorem}

This statement is sharp in the sense that for every $\theta\in
(0,\frac{\pi}{2})$ there exists a nonzero subharmonic function $u$
such that $M_r (u) = {\rm O} ((1-r)^{-\frac{\pi}{\theta}} )$ as $
r\to 1-,$ $\lim_{z \to \xi} u(z)=0$ for every $\xi\in {\mathbb T}
\setminus \{1 \},$ and $u=0$ in $\Omega_{\theta}(0)$; see
\cite[p. 283-286]{BeCo88} for more details. Theorem \ref{dal} can
be considered as the limit case of Theorem \ref{beco} with
$\theta=\frac{\pi}{2}.$

One may ask whether Theorems \ref{dal} and \ref{beco} can be
improved if $u$ is a {\it harmonic} function. Since the
function from \eqref{u_0} is harmonic, Theorem \ref{dal} is sharp
also in the class of harmonic functions. We note that Theorem
\ref{dal} for harmonic functions was  proved by V. L. Shapiro
\cite{Sh63} with a proof different from Dahlberg's proof.

On the other hand, Theorem \ref{beco} can be improved for harmonic
functions. The best results in this direction belong to F. Wolf;
see \cite{Wo41} and in particular Theorem $7.4.1$ and the proof of
Theorem $7.4.6$ therein.

\begin{theorem}[Wolf] \label{wolf}
Let $u$ be a harmonic function on the unit disc, and let $\theta \in
(0,\frac{\pi}{2})$ be fixed. Assume that
\begin{eqnarray*}
(i) & & \lim_{z\to e^{i\varphi}
\atop z \in \Omega_{\theta}(\varphi)} u(z) = 0 \mbox{ for every } \varphi \in [0,2\pi );\\
(ii) & &
M_r (u) = O ( e^{\varepsilon (1-r)^{-\frac{\pi}{2\theta}}} ) \quad
\text{ as } r\to 1-, \mbox{ for every }\, \varepsilon >0.
\end{eqnarray*}
Then $u=0$.
\end{theorem}

We modify an example of Rudin in \cite[p. 39]{AsBr91} to show that Theorem \ref{wolf} is optimal with
respect to the growth of $M_r(u)$, see Example \ref{app}.

Having in mind applications of uniqueness theorems to the study of
spaces of functions with integral norm, such as e.g. Hardy spaces
or Bergman spaces, one might wish to replace the $\sup$-norm in
the definition of $M_r (u)$ by the $L^p$-norm for some $p\geq 1$.
Unfortunately, such a generalisation appears to be either
impossible or superfluous. Indeed, as it was observed in
\cite{Sh63} for the function $u$ from \eqref{u_0}, $\|u (r
e^{i\cdot}) \|_p={\rm O}((1-r)^{-2+\frac{1}{p}}) = {\rm
o}((1-r)^{-2})$ as $r \to 1-.$ On the other hand, if $\|u (r
e^{i\cdot})\|_p = {\rm o}((1-r)^{-2+\frac{1}{p}})$ for some
subharmonic function $u$, then $\|u(r e^{i\cdot}) \|_{\infty} =
{\rm o}((1-r)^{-2})$ as $r \to 1-$ by simple estimates of Poisson
integrals, see \cite{Sh63}. So in this case, the conclusion $u=0$
follows from Theorem \ref{dal}. It is not clear what is  a natural
reformulation of Theorem \ref{beco} in the $L^p$-setting.

However, one can combine conditions from  Theorem \ref{beco} with
$L^p$-conditi\-ons by assuming that
$$
|u|\leq f \cdot g \text{ in } \D ,
$$
where $f$ is a lower semicontinuous function satisfying condition
(i) of Theorem \ref{beco}, while $g$ is measurable and $\sup_{r\in
[ 0,1)}\|g (r e^{i\cdot})\|_p < \infty.$ If, in addition, $M_r(u)$
is of at most polynomial growth, then we are able to prove
uniqueness theorems similar to the results stated above. This will
be the main subject of the present paper. For example, we prove
the following uniqueness theorem for harmonic functions.

\begin{theorem}[Uniqueness principle] \label{wmp}
Let $u$ be a harmonic function on the unit disc. Assume that
\begin{enumerate}
\item there exist a lower semicontinuous function $f : \D \to \R_+$ and a measurable function $g : \D\to\R_+$ such that
\begin{eqnarray*}
(i) & & |u(z)| \leq f(z) \, g(z) \quad \text{for every } z\in\D , \\[2mm]
(ii) & & \lim_{z\to e^{i\varphi} \atop z \in \Omega_{\theta}(\varphi)} f(z) = 0 \mbox{ for every } \varphi\in [0,2\pi ) \text{ and some } \theta \in (0,\frac{\pi}{2}) , \\
(iii) & & \sup_{r\in (0,1)} \int_0^{2\pi } g(r e^{i\varphi} ) \;
d\varphi <\infty ,
\end{eqnarray*}
\item there exists $m\geq 0$ such that
$$
M_r(u) = O ( (1-r)^{-m} ) , \quad r\to 1- .
$$
\end{enumerate}
Then $u=0$.
\end{theorem}

In this article we also obtain a version of Theorem \ref{wmp} for
non-negative subharmonic functions, as well as two versions of
these theorems which deal only with the {\em local} boundary
behaviour of harmonic and subharmonic functions. These main
theorems as well as a generalized Phragm\'en-Lindel\"of principle
are stated in Section \ref{sechf} and proved in Sections \ref{secpf1} and \ref{secpf2}.\\

In Section \ref{secapp} we present some applications of our main
results. First, we apply our local theorems to the
uniqueness theorem above, and to analytic continuation problems. After that, 
we obtain new optimal conditions for stability of operator
semigroups on Banach spaces with Fourier type. In
this application to operator semigroups conditions, using the  resolvent identity,
we obtain in a natural way conditions like those of 
Theorem~\ref{wmp}.

Finally, in Section \ref{appe} we sketch an example that shows how sharp is Theorem \ref{wolf}.\\

Note that uniqueness theorems for harmonic functions are also 
related to  the study of uniqueness problems for trigonometric
series. We do not discuss such relations here and refer the
interested reader to the papers \cite{Wo39, Wo41, Sh63, He71, AsBr91}.\\

The authors are grateful to N.~Nikolski for his attention to this work.

\section{Boundary behaviour of (sub-)harmonic functions} \label{sechf}

To formulate our main result, we define the rectangle
$$
R:=\{ z\in\C : -1 < {\rm Re}\, z  < 1, \, 0 < {\rm Im}\, z < 1 \} ,
$$
and for every $\theta\in (0,\frac{\pi}{2})$ we define the sector
$$
\Sigma_\theta := \{ z\in\C : \theta < {\rm arg}\, z < \pi -\theta \} .
$$
The statement below is a local halfplane version of the uniqueness
Theorem \ref{wmp}.

\begin{theorem} \label{main}
Let $u$ be a non-negative subharmonic function on $R$ which extends
continuously to $\bar{R} \setminus \R$. Assume that
\begin{enumerate}
\item \label{2} there exist a lower semicontinuous function $f : \bar{R}\setminus\R \to \R_+$ and a measurable function $g : R \to\R_+$ such that
\begin{eqnarray*}
(i) & & u(z) \leq f(z) \, g(z), \quad z\in R , \\[2mm]
(ii) & & \lim_{z\to \alpha \atop z\in \alpha+\Sigma_{\theta} } f(z) = 0 \mbox{ for every } \alpha\in (-1,1)
\mbox{ and some } \theta >0 , \\
(iii) & & \sup_{\beta \in (0,1)} \int_{-1}^1 g(\alpha + i\beta ) \; d\alpha<\infty ,
\end{eqnarray*}
\item \label{1} there exists $m\in [0,\frac{\pi}{\theta})$ such that
$$
\sup_{\alpha\in [-1,1 ]} u(\alpha +i\beta ) = O ( \beta^{-m} ) , \quad \beta\to 0+ .
$$
\end{enumerate}
 Then the function $u$ admits a continuous extension
to the interval $(-1,1)$, and $u=0$ on $(-1,1)$.
\end{theorem}

\begin{remark} Note that $\Sigma_\theta$ has an angular
opening $\pi-2\theta.$ Thus, Theorem \ref{main} agrees with
Theorem \ref{beco} as far as the relation between $\theta$ and $m$
is concerned. It is not clear whether in condition \eqref{1} one
can replace the ${\rm O}(\beta^{-m})$ for some $m\in
[0,\frac{\pi}{\theta})$ by just ${\rm o}(
\beta^{-\frac{\pi}{\theta}})$. By the argument from
\cite[p.283-286]{BeCo88} we cannot write ${\rm O}
(\beta^{-\frac{\pi}{\theta}})$, even for  $g=1$.
\end{remark}

For harmonic functions, the above result can be improved in the sense that in the second condition the constant $m$ does not depend on the angle $\theta$ from the first
condition.

\begin{theorem} \label{main2}
Let $u: R \to \C$ be a harmonic function which extends
continuously to $\bar{R} \setminus \R$. Assume that the condition
\eqref{2} of Theorem \ref{main} holds, and that, in addition,
\begin{enumerate} \setcounter{enumi}{1}
\item \label{harm1} there exists $m\geq 0$ such that
\begin{equation}\label{harmg}
\sup_{\alpha\in [-1,1 ]} | u(\alpha +i\beta ) | = O ( \beta^{-m} )
, \quad \beta\to 0+ .
\end{equation}
\end{enumerate}

Then the function $u$ admits a continuous extension to the interval $(-1,1)$, and $u=0$ on $(-1,1)$.
\end{theorem}

Clearly, the condition \eqref{harm1} of Theorem \ref{main2} is
weaker than the corresponding condition from Theorem \ref{main},
so that we obtain a stronger result for harmonic functions. The
proof of Theorem \ref{main2} is similar to the proof of Theorem
\ref{main} apart from one step. In this step, we need a
proposition of Phragm\'en-Lindel\"of type which we state
separately.

\begin{proposition} \label{prop}
Let $u: \Delta_{\rho} \to \R$ be a harmonic function on the finite
open sector
$$
\Delta_{\rho} := \{ z\in\C : 0<\arg z <\frac{\pi}{2} , \, |z|<
\rho\} .
$$
Assume that
\begin{enumerate}
\item \label{a1} $u$ has a continuous extension to $\overline{\Delta_{\rho}}\setminus \{ 0\}$,
\item \label{a2} there exist $l\in [0,2)$, $\theta\in (0,\frac{\pi}{2})$ such that $\frac{\theta}{\pi}\not\in\Q$ and
$$
\sup_{z\in \overline{\Delta_{\rho}}\setminus \{ 0\} \atop \arg
z\in \{ 0 , \, \theta , \, \frac{\pi}{2}\} } |z|^l\, |u(z)| <
\infty ,
$$
\item \label{a4} there exists $m\in [0,\infty )$ such that
$$
\sup_{z\in \Delta_{\rho}} |{\rm Im}\, z|^m\, |u(z)| <\infty .
$$
\end{enumerate}
Then
$$
\sup_{z\in \Delta_{\rho}} |z|^l\, |u(z)| <\infty .
$$
\end{proposition}

\begin{remark}
Classical Phragm\'en-Lindel\"of principle
does not directly apply to the situation of Theorem \ref{prop} since
there is no relation between the numbers $\theta$ and $m$ in
the assumptions \eqref{a2} and \eqref{a4}. 
\end{remark}

\section{Proof of Theorems \ref{main} and \ref{main2}} \label{secpf1}

In what follows, we denote by $C$ a positive constant which may vary from line to line.\\

We start by stating two lemmas which are needed for the proofs of Theorems \ref{main}, \ref{main2} and Proposition \ref{prop}.\\

The first lemma, a version of the mean value inequality for
subharmonic functions, can be found in {\sc Koosis} \cite[Lemma,
Section VIII D.2]{Ko80} or {\sc Garnett} \cite[Chapter III, Lemma
3.7]{Ga81}.

\begin{lemma} \label{lmean}
Let $u$ be a non-negative subharmonic function defined on a domain
$D\subset \C$ and let $p > 0.$ Then for every $z \in D$ and every
closed ball $\bar{B}(z,r)\subset D$ we have
\begin{equation}\label{mean}
u(z)^p  \le C_p \frac{1} {|B(z,r)|} \int_{B(z,r)} u(y)^p  \; dy ,
\end{equation}
where $C_p$ is a constant depending only on $p.$
\end{lemma}

The second lemma is a variant of a result by {\sc Domar} \cite{Do58}. Our proof uses an idea from \cite[Lemma 4]{AsBr91}. 

\begin{lemma} \label{domar2}
Let $u: U \setminus \{ 0\} \to \R_+$ be a subharmonic function,
where $U$ is a neighbourhood of the closed unit disc $\bar{\D}$.
Assume that for some constants $C$, $m\ge 0$ and for all $z\in
\D\setminus \{ 0\}$
\begin{equation} \label{a41}
u(z) \leq C\, |\Im z |^{-m} .
\end{equation}
Then there is another constant $C>0$ such that for
all $z\in \D\setminus \{ 0\}$
$$
u(z) \leq C \, |z|^{-m} .
$$
\end{lemma}

\begin{proof}
Let $z\in \D$ be fixed, and let $r := r(z) := |z|$.

If $|\Im z| \geq r/2$, then the above inequality is a direct consequence of assumption \eqref{a41}.

So we can assume that $|\Im z| \leq r/2$. Choose $\delta > 0$ independent of $z\in \D\setminus \{ 0\}$ such that the closed disc $\bar{B} (z, \delta |z|)$ is contained in $U\setminus \{ 0\}$. Applying the mean value inequality \eqref{mean} with $p=\frac{1}{2m}$ (Lemma \ref{lmean}) and using the assumption \eqref{a41} we obtain
\begin{eqnarray*}
u(z)^{\frac{1}{2m}} & \leq & C \frac{1}{\delta^2 r^2 \pi} \int_{B(z, \delta r)}
u(y)^{\frac{1}{2m}} \; dy \\
& \leq & \frac{C}{r^2} \int_{B(z, \delta r)} |{\rm Im}\, y|^{\frac{1}{2}} \; dy \\
& \leq & \frac{C}{r^2} \int_{\Re z - \delta r}^{\Re z +\delta r}
\int_{\Im z - \delta r}^{\Im z + \delta r} |\beta|^{-\frac{1}{2}}
\;
d\beta \; d\alpha \\
& \leq & \frac{C}{r} \left( |\Im z + \delta r|^{\frac{1}{2}} + |\Im z - \delta r|^{\frac{1}{2}} \right) \\
& \leq & C \, r^{-\frac{1}{2}} \quad = \quad C \, |z|^{-\frac12} .
\end{eqnarray*}
This is the claim.
\end{proof}

For the proof of Theorem \ref{main} we also need the following
Phragm\'en-Lindel\"of principle for functions which are
subharmonic in a finite sector. Similar (and more general)
statements for infinite sectors can be found in {\sc Levin}
\cite[Theorem 3, p. 49]{Le96} or {\sc Hardy \& Rogosinski}
\cite{HaRo46}. The proof of our statement
is an easy adaptation of their proofs.

\begin{proposition}[Classical Phragm\'en-Lindel\"of principle] \label{pl}
Let $u$ be a  non-negative subharmonic function on the open finite
sector
$$
\Delta_{\theta ,r} := \{ z\in\C : |\arg z |<\theta , \, |z|< r
\} , \quad \theta \in (0,\pi ), \, r>0.
$$
Assume that
\begin{enumerate}
\item \label{p1} $u$ has a lower semicontinuous extension to $\overline{\Delta_{\theta ,r}}\setminus \{ 0\}$,
\item \label{p2} there exist $0\leq l\leq m < \frac{\pi}{2\theta}$ such that
$$
\sup_{z\in \Delta_{\theta ,r} \atop |\arg z| = \theta } |z|^l\,
u(z) < \infty ,
$$
 and
$$
\sup_{z\in \Delta_{\theta ,r}} |z|^m\, u(z) <\infty .
$$
\end{enumerate}
Then
$$
\sup_{z\in \Delta_{\theta ,r}} |z|^l\, u(z) <\infty .
$$
\end{proposition}

\begin{proof}[Proof of Theorem \ref{main}]
We call a point $\alpha\in (-1,1)$ {\em regular} if $u$ has a
continuous extension to a neighbourhood of $\alpha$ in $\bar{R}$
and we call it {\em singular} otherwise. We denote the set of all
singular points in $(-1,1)$ by $S$. Clearly, $S$ is closed in
$(-1,1)$.

{\bf Step 1:} We first prove that $u=0$ on $(-1,1)\setminus S$.
Otherwise, by continuity, we can find a
nonempty interval $U\subset (-1,1)\setminus S$ such that $u(\alpha
) \geq c >0$ for every $\alpha\in U$. By assumption \eqref{2} (ii)
and Egorov's theorem, we find a set $V\subset U$ of positive
measure such that
$$
\lim_{\beta\to 0+} \sup_{\alpha\in V} f(\alpha +i\beta ) = 0 .
$$
Together with assumptions \eqref{2} (i), (iii) and Fatou's lemma this implies
\begin{eqnarray*}
0 & < & \liminf_{\beta\to 0+} \int_V u(\alpha +i\beta ) \; d\alpha \\
& \leq & \limsup_{\beta\to 0+} \int_V f(\alpha +i\beta ) \, g(\alpha +i\beta ) \; d\alpha \quad = \quad 0 ,
\end{eqnarray*}
a contradiction. Hence, $u$ vanishes at every regular point and it remains to prove that the set $S$ of singular points in $(-1,1)$ is empty.\\

We assume in the following that $S$ is nonempty and we will show that this leads to a contradiction.\\

{\bf Step 2:} We define for every $n\in\N$
\begin{eqnarray*}
S_n & := & \{ \alpha\in S : \sup_{z\in R\cap \alpha +\Sigma_{\theta}} f(z) \leq n \} \\
& = & \bigcap_{\beta >0} \{ \alpha\in S : \sup_{z\in R\cap \alpha +\Sigma_{\theta} \atop {\rm Im}\, z\geq \beta} f(z) \leq n \} .
\end{eqnarray*}
By lower semicontinuity of the function $f$, the sets $S_n$ are closed. By assumption \eqref{2} (ii),
$$
S = \bigcup_{n\in\N} S_n .
$$
By Baire's category theorem, there exists $n\in\N$ such that
$S_{n}$ has nonempty interior in $S$, i.e. there exists a nonempty
interval $(a,b)\subset (-1,1)$ such that
\begin{equation} \label{e2}
S\cap (a,b) = S_{n} \cap (a,b) \not= \emptyset .
\end{equation}
Without loss of generality we can assume that $a> -1$ and $b< 1$. Since the set
$(a,b)\setminus S$ is open, it is the countable union of mutually
disjoint intervals $I_k:= (a_k ,b_k)$, 
$$
(a,b)\backslash S= \bigcup_k \; (a_k,b_k),
$$

In the following, we put for every $k$
\begin{eqnarray*}
l_k & := & b_k -a_k = |I_k| \qquad \text{ and} \\
h_k & := & l_k \, \frac{\tan \theta}{2} ,
\end{eqnarray*}
and we also define
$$
h := \min \{ (a+1) \tan \theta , \, (1-b) \tan\theta , \, 1 \} > 0 .
$$
Furthermore, we consider the rectangles
$$
\tilde{R} := \{ z\in R : a < {\rm Re}\, z < b \}
$$
and
$$
R_k := \{ z\in R : a_k <{\rm Re}\, z<b_k \} .
$$
The parameter $h$ is only of technical interest: it will be convenient to know that for every $\alpha\in [a,b]$ one has
$$
\{ z\in\C : z\in \alpha +\Sigma_\theta , \, {\rm Im}\, z < h \} \subset R .
$$

{\bf Step 3:} We prove that
\begin{equation} \label{e3}
\sup_{0< \beta < 1} \int_a^b u(\alpha +i\beta )\; d\alpha <\infty .
\end{equation}

{\bf Step 3.1:} It follows from \eqref{e2} that for every $\beta\in (0,1)$
$$
\int_{(a,b)\cap S} u(\alpha +i\beta )\; d\alpha \leq n \, \int_{(a,b)\cap S} g(\alpha +i\beta )\; d\alpha \leq n\, C ,
$$
so that
\begin{equation} \label{e31}
\sup_{0< \beta < 1} \int_{(a,b)\cap S} u(\alpha +i\beta ) \; d\alpha < \infty .
\end{equation}

{\bf Step 3.2:} Fix $\beta\in (0,1)$. If $\beta \geq h_k$, then
$$
I_k + i\beta \subset (S\cap (a,b)) + \Sigma_{\theta} ,
$$
so that
\begin{eqnarray*}
\sum_{k \atop \beta \geq h_k} \int_{a_k}^{b_k} u(\alpha +i\beta ) \; d\alpha & \leq & \sum_{k \atop \beta \geq h_k} \int_{a_k}^{b_k} n \, g(\alpha +i\beta ) \; d\alpha \\
& \leq & n \, \sup_{0< \beta <1} \int_a^b g(\alpha +i\beta ) \; d\alpha \quad \leq \quad n C .
\end{eqnarray*}
Hence,
\begin{equation} \label{e32}
\sup_{0<\beta< 1} \sum_{k \atop \beta \geq h_k} \int_{a_k}^{b_k} u(\alpha +i\beta )\; d\alpha <\infty .
\end{equation}

{\bf Step 3.3:} We prove that
\begin{equation} \label{e33}
\sup_{0<\beta< 1} \sum_{k \atop \beta < h_k} \int_{a_k}^{b_k} u(\alpha +i\beta ) \; d\alpha <\infty .
\end{equation}

In what follows, if $D\subset\C$, $D\not=\C$, is a simply
connected domain then we denote by $\omega_{D}$ the harmonic
measure associated with $D$; see e.g. \cite{Ko88}. If $h:\partial
D \to \R_+$ is a Borel function such that $\int_{\partial D}
|h(\xi )|\; d\omega_D (\xi ,z) <\infty$ for every $z\in D$, then
we can define
\begin{equation}\label{poi1}
P_{D} (h) (z) := \int_{\partial D} h(\xi ) \; d\omega_{D} (\xi ,z) , \quad z\in D .
\end{equation}
The function $P_D (h)$ is harmonic in $D$, and if $h$ is continuous and bounded, then
$P_D (h)$ is actually the unique continuous solution of the Dirichlet problem on $\bar{D}$ with boundary value $h$.\\
If $\partial D$ is piecewise smooth, and  $\varphi: D \mapsto \D$
is a conformal mapping, then \eqref{poi1}
can be rewritten as
\begin{equation}\label{poi2}
P_D (h) (\varphi^{-1}(z)) =\frac {1}{2\pi} \int_{\partial \D} h(\varphi^{-1}(\xi))
\frac{1-|z|^2}{|z-\xi|^2}\, |d\xi|.
\end{equation}

The next technical result  will be used several times in the sequel.

\begin{lemma}\label{harm}
Let $D\subset\C$ be a piecewise smooth, bounded, simply connected domain such
that $\partial D$ has a finite set of corners $C.$ Assume that the
opening of each corner is $\pi\alpha$, $\alpha \in [\frac{1}{2},
1).$ Let $0<\alpha\beta<1$, and let $h:\partial D \to \mathbb C$ be a function which is continuous on
$\partial D\setminus C$ and which satisfies
$$
|h(\xi)|={\rm O} \left(|\xi-c|^{-\beta} \right),\qquad \text{as } \xi \to c, \, \xi \in \partial D,
$$ 
for each $c \in C.$ Then the following statements are true:
\begin{itemize} 
\item [(i)] for every $z \in D$ one has $h \in L^1 (\partial D, d\omega_D (\cdot, z))$;
\item [(ii)] if $h$ is continuous at $\xi \in \partial D,$ then $\lim_{z \to \xi} P_D(h)(z)=h(\xi);$
\item [(iii)] $|P_D (h)(z)|= {\rm O} \left(|z-\xi|^{-\beta}
\right)$ as $z \to \xi$, $\xi \in C.$
\end{itemize}
\end{lemma}

Lemma \ref{harm} follows from \eqref{poi2} and the facts that
$\varphi$ extends continuously to $\bar D$, satisfies
$\varphi(\xi)={\rm O} (|\xi-c|^{1/\alpha})$ as $\xi\to c,$ for every $c \in C$,
and $\varphi^{-1}$ extends continuously to $\bar \D$.\\

{\bf Step 3.3.1:} We show that for every $k$
\begin{equation} \label{e330}
u(z) \leq C(k)\, \max \big\{ \frac{1}{|z-a_k|} , \frac{1}{|z-b_k|} \big\} \text{ for every } z\in R_k .
\end{equation}

Fix $k$, let $\theta_1 \in (\theta ,\frac{\pi}{2} )$ and choose $\delta>0$ so small that
$$
z\in \Sigma_{\theta_1} \quad \Rightarrow \quad B(z,|z|\delta ) \subset \Sigma_{\theta} .
$$
Then for every $z\in \Sigma_{\theta_1}$ with $|z|$ small enough the mean value inequality implies
\begin{eqnarray*}
u(a_k+z) & \leq & \frac{C}{|z|^2} \int_{B(a_k+z,|z|\delta )} u(z') \; dz' \\
& \leq & \frac{Cn}{|z|^2} \int_{B(a_k+z,|z|\delta )} g(z') \; dz' \\
&\leq & \frac{Cn}{|z|^2} \int_{{\rm Im}\, z - |z| \delta}^{{\rm Im}\, z + |z|\delta} \int_{a}^{b} g(\alpha +i\beta) \, d\alpha \, d\beta \\
 & \leq & \frac{C}{|z|} ,
\end{eqnarray*}
for some constant $C$ which depends only on $\delta$ (i.e. $\theta_1$ and $\theta$), $n$ and $g$.
Similarly, for every $z\in \Sigma_{\theta_1}$ with $|z|$ small enough,
$$
u(b_k +z) \leq \frac{C}{|z|}.
$$
Now we choose $\theta_1$ sufficiently close to $\theta$ so that $m<\frac{\pi}{\theta_1}$
(where $m$ is as in assumption \eqref{1}). The Phragm\'en-Lindel\"of principle from Proposition \ref{pl} implies \eqref{e330}.\\

{\bf Step 3.3.2:} By \eqref{e330} and Lemma \ref{harm}, for every
$k$ and every $z\in R_k$ the function $u|_{\partial R_k}$ is
integrable with respect to the harmonic measure $\omega_{R_k}
(\cdot ,z)$. We show that for every $k$
\begin{equation} \label{e331}
u(z) \leq P_{R_k} (u) (z) \text{ for every } z\in R_k .
\end{equation}
Clearly, the function $v_k (z) := u(z) - P_{R_k} (u)(z)$ ($z\in
R_k$) is subharmonic in $R_k$, continuous up to
$\overline{R_k}\setminus\{ a_k ,b_k\}$, and $v_k = 0$ on $\partial
R_k\setminus \{ a_k ,b_k\}$. Let $R_{k,\varepsilon} := \{ z\in R_k
: {\rm Im}\, z>\varepsilon \}$. Then, for every $z\in R_k$ and
every $\varepsilon\in (0,{\rm Im}\, z)$,
$$
v_k (z) \leq \int_{\partial R_{k,\varepsilon}} v_k (\xi ) \; d\omega_{R_{k,\varepsilon}} (\xi ,z) .
$$
Letting $\varepsilon\to 0+$, by \eqref{e330} and the dominated convergence theorem, we obtain
$$
v_k (z) \leq \int_{\partial R_{k}} v_k (\xi ) \; d\omega_{R_{k}} (\xi ,z) = 0 .
$$
This implies \eqref{e331}.

{\bf Step 3.3.3:} By \eqref{e331},
\begin{eqnarray}
\nonumber \sum_{k \atop \beta < h_k} \int_{a_k}^{b_k} u(\alpha +i\beta ) \; d\alpha & \leq & \sum_{k \atop \beta < h_k} \left\{ \int_{a_k}^{b_k} P_{R_k} (u\, \chi_{L_k} ) (\alpha +i\beta ) \; d\alpha + \right. \\
\label{e332} & & \phantom{+++} + \int_{a_k}^{b_k} P_{R_k} (u \, \chi_{M_k} ) (\alpha +i\beta ) \; d\alpha + \\
\nonumber & & \phantom{+++} \left. + \int_{a_k}^{b_k} P_{R_k} (u\,
\chi_{M_k'} ) (\alpha +i\beta ) \; d\alpha \right\} ,
\end{eqnarray}
where
\begin{eqnarray*}
L_k & := & \{ z\in \partial R_k : {\rm Im}\, z =1 \} , \\
M_k & := & \{ z\in\partial R_k : {\rm Re}\, z = a_k \} , \text{ and} \\
M_k' & := & \{ z\in\partial {R_k} : {\rm Re}\, z = b_k \} .
\end{eqnarray*}

By continuity, we have
$$
\sup_{\alpha\in [-1,1]} u(\alpha +i) = C <\infty,
$$
and
$$
P_{R_k} (u\, \chi_{L_k}) (z) \leq C \text{ for every } z\in R_k \text{ and every } k .
$$
Hence,
$$
\sum_{k \atop \beta < h_k} \int_{a_k}^{b_k} P_{R_k} (u\, \chi_{L_k} ) (\alpha +i\beta ) \; d\alpha \leq \sum_{k \atop \beta < h_k} \int_{a_k}^{b_k} C \; d\alpha \leq C\, (b-a) ,
$$
so that we have estimated the first term on the right-hand side of
\eqref{e332} by a constant independent of $\beta\in (0,h_k)$.

Let
$$
Q := \left \{ z\in\C : 0 < \arg z < \frac{\pi}{2} \right \} .
$$
For every $k$ we define $u_k : \partial Q \to \R_+$ by
$$
u_k (i\beta ) := \left\{ \begin{array}{ll}
u (a_k +i\beta )  & \text{if } 0 <\beta \leq 1 , \\[2mm]
0 & \text{if } \beta >1 ,
\end{array} \right. \text{ and } u_k (\alpha ) = 0 , \quad \alpha >0 .
$$

Fix $\beta > 0$ and put $r = \beta \, \cos \theta$. Then one
checks that the intervals $((a_k-r,a_k+r))_{\!\!\!\!\! k \atop \!\!\!\beta <
h_k}$ are mutually disjoint. This fact, the mean value inequality, and
the assumption \eqref{2} imply that
\begin{eqnarray*}
\sum_{k \atop \beta < h_k} u_k (i\beta ) & = & \sum_{k \atop \beta < h_k} u (a_k +i\beta ) \\
& \leq & \frac{C}{r^2} \sum_{k \atop \beta < h_k} \int_{|z-a_k-i\beta|\leq r} u(z) \; dz \\
& \leq & \frac{Cn}{r^2} \sum_{k \atop \beta < h_k} \int_{\beta -r}^{\beta +r} \int_{a_k -r}^{a_k+r} g(\alpha'+i\beta') \; d\alpha' \; d\beta' \\
& \leq & \frac{C}{r^2}  \int_{\beta -r}^{\beta +r} \int_a^b g(\alpha'+i\beta' ) \; d\alpha' \; d\beta' \\
& \leq & \frac{C}{r} \quad \leq \quad \frac{C}{\beta} .
\end{eqnarray*}

We deduce that
\begin{eqnarray*}
\sum_{k \atop \beta < h_k} \int_{a_k}^{b_k} P_{R_k} (u\, \chi_{S_k} ) (\alpha +i\beta ) \; d\alpha & \leq & \sum_{k \atop \beta < h_k} \int_0^\infty P_Q (u_k ) (\alpha +i\beta ) \; d\alpha \\
& = & \int_0^\infty P_Q (\sum_{k \atop \beta < h_k} u_k ) (\alpha +i\beta ) \; d\alpha \\
& \leq & C\, \int_0^\infty P_Q (v ) (\alpha +i\beta ) \; d\alpha ,
\end{eqnarray*}
where $v:\partial Q\to\R_+$ is defined by
$$
v(i\beta' ) := \frac{1}{\beta'} \text{ and } v(\alpha' ) := 0 , \quad \alpha', \, \beta' >0 .
$$
Note that
$$
P_Q (v) (z) = P(z) , \quad z\in Q ,
$$
where $P$ is the Poisson kernel in the upper half-plane $\mathbb
C_+ := \{ z\in\C : {\rm Im}\, z >0\}$. Indeed, if $v_0 =|P_Q(v) - P|,$ then $v_0$ is subharmonic,
continuous in $\overline Q \setminus \{0 \},$ and, by Lemma \ref{harm} (iii), $|v_0(z)|={\rm O} (|z|^{-1})$ as $z \to 0.$ By Proposition \ref{pl}, $v_0$ is bounded, and by the classical maximum principle we conclude that $v_0=0.$

Since
$$
\int_\R P(\alpha +i\beta ) \; d\alpha = 1 \text{ for every } \beta >0 ,
$$
we have proved that the second term on the right-hand side of
\eqref{e332} can be estimated by a constant independent of
$\beta\in (0,1)$. The same argument works for the third term on
the right-hand side of \eqref{e332}.

Hence, the right-hand side of \eqref{e332} is uniformly bounded in $\beta\in (0,1)$, and this proves \eqref{e33}.

Summing up \eqref{e31}, \eqref{e32} and \eqref{e33}, we obtain \eqref{e3}.

{\bf Step 4:} By \eqref{e3} and the weak$^*$ sequential
compactness of the unit ball of $M ([a,b]) = C([a,b])^*$, there
exists a sequence $(\beta_k)\searrow 0, k \to \infty,$ and a
finite positive Borel measure $\mu\in M ([a,b])$ such that
$$
\underset{k\to\infty}{w^*-\lim}\,\, u(\cdot +i\beta_k ) = \mu \text{ in } M
([a,b]) .
$$

Let $P_{\C_+} (\mu )$ be the Poisson integral of $\mu$ in the upper half-plane $\C_+$.

Let $L:=\partial{\tilde{R}}\cap \mathbb C_+$ be the part of the
boundary of $\tilde{R}$ which does not lie on the real axis. Since
$$
u(z) \leq C\, \max \left \{ \frac{1}{|z-a|} , \frac{1}{|z-b|}
\right \} , \quad z\in L ,
$$
by Lemma \ref{harm}, the function $P_{\tilde{R}} (u\, \chi_L )$ is
well defined and harmonic in $\tilde{R}$. Moreover,
$$
\lim_{\beta\to 0+} P_{\tilde R} (u\, \chi_L ) (\cdot +i\beta ) = 0 \text{ locally uniformly on } (a,b) .
$$
Hence, if we define $v:= \max \left((u-P_{\tilde{R}} (u\, \chi_L
)), 0\right)$, then $v$ is a non-negative subharmonic function on
$\tilde{R}$, continuous up to $L$, and $v=0$ on $L$. Moreover, for
every continuous $\varphi$ with support in $(a,b),$
\begin{equation}\label{A}
\lim_{k\to\infty} \int_a^b v(\alpha +i\beta_k ) \varphi (\alpha )
\; d\alpha = \int_a^b \varphi (\alpha ) \; d\mu .
\end{equation}
On the other hand, for every such $\varphi$ we also have
\begin{equation}\label{B}
\lim_{\beta\to 0+} \int_a^b P_{\mathbb C_+} (\mu ) (\alpha +i\beta
) \varphi (\alpha ) \; d\alpha = \int_a^b \varphi (\alpha ) \;
d\mu.
\end{equation}
Now, given $\beta>0$, denote $\tilde R_{\beta}:=i\beta+\tilde R$. For $z\in\tilde R$,
\begin{gather*}
P_{\mathbb C_+}(\mu)(z)=\lim_{\beta \to 0+}P_{\mathbb C_+}(\mu)(z+i\beta)\\
\ge \limsup_{\beta \to 0+}\int_{a}^{b} P_{\mathbb C_+}(\mu)(\alpha + i\beta) \, 
d\omega_{\tilde R_{\beta}}(\alpha,z+i\beta)\\=
\limsup_{\beta \to 0+}\int_{a}^{b} P_{\mathbb C_+}(\mu)(\alpha + i\beta) \, d\omega_{\tilde R} (\alpha, z). 
\end{gather*}
 
Furthermore, since the function
$$
\alpha \mapsto \frac{d \omega_{ \tilde R} (\alpha, z)}{d \alpha}
$$
is continuous on $[a,b]$ and vanishes at $a$ and $b$, and
since
$$
\int_{a}^{b}P_{\mathbb C_+}(\mu )(\alpha + i\beta)\,d\alpha\le C \text{ for every } \beta >0,
$$
we obtain by \eqref{B} that
$$
P_{\mathbb C_+} (\mu ) (z) \ge \limsup_{\beta \to 0+}\int_{a}^{b} \frac{d\omega_{\tilde R} (\alpha, z )}{d \, \alpha} \, d \mu(\alpha ).
$$
On the other hand, if we set $R^*_{\beta_k}:=\{z \in \tilde R: {\rm Im}\, z > \beta_k \}$, then
\begin{gather*}
v(z)=\lim_{k\to\infty}v(z+i\beta_k)\le \liminf_{k \to \infty} 
\int_{a}^{b} v(\alpha + i\beta_k)\, d\omega_{R^*_{\beta_k}} (\alpha,z+\beta_k)\\ \le
\liminf_{k \to \infty} 
\int_{a}^{b} v(\alpha + i\beta_k)\, d\omega_{\tilde R_{\beta_k}} (\alpha,z+\beta_k)\\=
\liminf_{k \to \infty} \int_{a}^{b} v(\alpha + i\beta_k)\, d\omega_{\tilde R} (\alpha, z).
\end{gather*}
 
As above, we conclude by \eqref{e3} and by \eqref{A} that
$$
v(z) \le \liminf_{k \to \infty} \int_{a}^{b} \frac{d\omega_{\tilde
R} (\alpha, z)}{d \alpha} \, d\mu(\alpha).
$$
Thus,
\begin{equation}
P_{\C_+} (\mu ) (z) \geq v(z) \text{ for every } z\in\tilde{R} ,
\label{lbl}
\end{equation}
and it follows from \eqref{A} and \eqref{B} that
 for every $a<a'<b'<b$
\begin{equation} \label{e410}
\lim_{k\to\infty} \int_{a'}^{b'} | (P_{\mathbb C_+} (\mu ) - v)
(\alpha +i\beta_k )| \; d\alpha = 0 .
\end{equation}

Assume that $\mu ((a,b))\not= 0$. Then one finds $a<a'<b'<b$ such that $\mu ((a',b')) >0$. By assumption \eqref{2} (ii) and Egorov's theorem, there exists a set $T\subset (a',b')$ of positive $\mu$-measure such that
$$
\lim_{\beta\to 0+} \sup_{\alpha\in (-1,1)\atop \alpha +i\beta \in
T +\Sigma_{\theta}} f(\alpha +i\beta) = 0 .
$$
Together with assumption \eqref{2} (i) and (iii) this implies
$$
\limsup_{\beta\to 0+} \int_{\alpha\in (-1,1)\atop \alpha +i\beta
\in T +\Sigma_{\theta}} u(\alpha +i\beta )\; d\alpha = 0 ,
$$
and thus also
\begin{equation} \label{e411}
\limsup_{\beta\to 0+} \int_{\alpha\in (-1,1)\atop \alpha +i\beta
\in T +\Sigma_{\theta}} v(\alpha +i\beta )\; d\alpha = 0 .
\end{equation}
On the other hand, one easily checks that
$$
\int_{\alpha\in (-1,1)\atop \alpha +i\beta \in T +\Sigma_{\theta}}
P_{\mathbb C_+} (\mu) (\alpha +i\beta )\,d\alpha \ge C(\theta) \mu (T)> 0
$$
for every $\beta >0$. This inequality contradicts \eqref{e410} and \eqref{e411}, and therefore the assumption $\mu ((a,b)) \not= 0$ is wrong.

However, if $\mu ((a,b)) = 0$, then by \eqref{lbl}, $v=0$ on $\tilde R$, and $u$ is continuous on
$\tilde R\cup (a,b)$. Thus,
$$
S \cap (a,b) =\emptyset ,
$$
which contradicts to  \eqref{e2}.

Hence, the assumption that $S$ is nonempty is  wrong.
Therefore, $S$ is empty and the theorem is proved.
\end{proof}

\begin{proof}[Proof of Theorem \ref{main2}]
The proof of Theorem \ref{main} goes through for the non-negative
subharmonic function $|u|$, except for the inequality
\begin{equation} \label{e331a}
|u(z)| \leq P_{R_k} (|u|) (z) \text{\, for every } z\in R_k 
\end{equation}
from Step $3.3.2$ that needs a justification.

As in the Step 3.3.1 of the proof of Theorem \ref{main} one proves that for
every $\theta_1\in (\theta ,\frac{\pi}{2} )$ there exists a constant $
C\geq 0$ such that for every $z\in \Sigma_{\theta_1}$ with $|z|$ small enough one has
$$
|u(a_k +z)|\leq \frac{C}{|z|} \text{\, and \,} |u(b_k +z)|\leq \frac{C}{|z|} .
$$

By Proposition  \ref{prop}, there exists a constant $C=C(k)\geq 0$
such that for every $z\in R_k$,
$$
|u(z)| \leq C\, \max \left \{ \frac{1}{|z-a_k|} , \,
\frac{1}{|z-b_k|} \right \} .
$$
Arguing as in Step $3.3.2,$ we obtain \eqref{e331a} for harmonic
functions satisfying the weaker growth condition \eqref{harmg}. The
rest of the proof is the same as in Theorem \ref{main}.
\end{proof}

\section{Proof of Proposition \ref{prop}} \label{secpf2}

\begin{proof}[Proof of Proposition \ref{prop}]
Without loss of generality, we may in the following assume that $m$ is an integer larger than $l$.

{\bf Step 1:} By assumptions \eqref{a1}-\eqref{a2} and Lemma \ref{harm} (i),
the function $u|_{\partial \Delta_{\rho}}$  is integrable with
respect to the harmonic measure $\omega_{\Delta{_\rho}}$
associated with the sector $\Delta_{\rho}$. Therefore, the function
$$
v:= P_{\Delta_{\rho}} (u|_{\partial \Delta_{\rho}}
)=\int_{\partial \Delta_{\rho}} u(\xi ) \; d\omega_{\Delta_{\rho}}
(\xi ,z ) , \quad z\in \Delta_{\rho},
$$
is well-defined. Moreover, by the assumption \eqref{a1} and Lemma
\ref{harm} (ii), $v$ is continuous up to
$\overline{\Delta_{\rho}} \setminus\{ 0\}$ and, by the assumption
\eqref{a2} and Lemma \ref{harm} (iii), there exists a constant
$C\geq 0$ such that
\begin{equation} \label{vestimate}
|v(z)| \leq C |z|^{-l} \text{ for all } z\in \Delta_{\rho} .
\end{equation}

Hence, it suffices to show that $u=v$.

{\bf Step 2:} Let $w:= u-v$. We define the punctured disc
$$
D_\rho := \{ z\in\C : |z|< \rho \} \setminus \{ 0\} ,
$$
and we denote also by $w$ the harmonic function on $D_\rho$ which
one obtains by reflecting the original $w$ first at the real axis
(using the Schwartz reflection principle) and then at the
imaginary axis (using the Schwartz reflection principle again).

We prove that there exists a constant $C >0$ such that
\begin{equation} \label{domar}
|w(z) | \leq C \, |z|^{-m}  \quad \mbox{for all } z\in D_\rho .
\end{equation}
In fact, by assumption \eqref{a4} and the estimate \eqref{vestimate}, there exists a constant $C\geq 0$ such that
$$
|w(z)| \leq C \, |{\rm Im}\, z|^{-m} , \quad z\in R .
$$
The estimate \eqref{domar} follows then from a simple application of Lemma \ref{domar2}.

{\bf Step 3:} Since $0$ is an isolated singularity of the harmonic function $w$, we have
$$
w(r e^{i \phi})= a\,\log r+\sum_{n \in \Z}(r^n a_n \cos n \phi + r^n b_n\sin n\phi) 
$$
for all $z=re^{i\varphi} \in D_\rho$. The series converges absolutely on every compact subset of 
$D_\rho$. The estimate \eqref{domar} implies that
\begin{equation}\label{repre}
w(r e^{i \phi}) = a \, \log r + \sum_{-m \le n \le -1}( a_n
r^n\cos n \phi + b_n r^n\sin n \phi)+ w_0 (r e^{i \phi}),
\end{equation}
with $w_0$ bounded in $D_\rho$.

By assumption \eqref{a2}, the function $r \to r^lw(re^{i\theta})$ is bounded on $(0,1),$ and, 
moreover, $\cos n\theta$ and $ \sin n \theta, -m \le n \le -1,$ are nonzero. This implies that 
$a_n=b_n=0,-m\le n\le-2$. Now recall that $w=0$ on $(0,\rho)$ and on $(0,i\rho)$. Therefore, 
letting $\phi=0$ in the representation \eqref{repre}, we conclude that $a=0$ and $a_{-1}=0$. Finally, if we let $\phi = \frac{\pi}{2}$ in \eqref{repre}, we obtain that $b_{-1}=0$.

Hence, $w=w_0$ is bounded. The classical maximum principle (note that $w=0$ on $\partial{{D_r}}$) implies that $w=0$, i.e. $u=v$.
\end{proof}

\section{Applications}\label{secapp}

\subsection{Uniqueness principle}

As a first application we prove the uni\-queness principle for harmonic functions formulated
in the Introduction.

\begin{proof}[Proof of Theorem \ref{wmp}]
Let
$$
\varphi_\theta(z)=\exp[2\pi i(z+\theta)], \qquad z\in\C,\, 0\le\theta<1.
$$
This holomorphic function $\varphi_\theta$ maps the interior of the rectangle
$$
R:=\{ z\in\C : -1 < {\rm Re}\, z  < 1, \, 0 < {\rm Im}\, z < 1 \}
$$
onto the set
$$
\{z\in\C :e^{-2\pi}<|z|<1\}\setminus (0,1)\exp[2\pi \theta i].
$$
Moreover, $\varphi_\theta$ is everywhere locally invertible. The function 
$\tilde{u}_\theta=u\circ\varphi_\theta$ is harmonic and it satisfies all the assumptions from 
Theorem \ref{main2}. By Theorem \ref{main2}, $\tilde{u}_\theta$ has a continuous extension 
to the interval $(-1,1)$ and this extension is $0$ on that interval. By local invertibility of $\varphi_\theta$ this implies that the function $u$ has a continuous extension to the closed unit disc $\bar{\D}$ and that $u|_{\partial \D} =0$. The claim follows from the classical maximum principle.
\end{proof}

A similar uniqueness principle holds for non-negative subharmonic functions.
The proof is similar to that of Theorem \ref{wmp}, and is based on Theorem \ref{main}.

\begin{theorem}[Uniqueness principle] \label{wmp2}
Let $u$ be a non-negative subharmonic function on the unit disc. Assume that
\begin{enumerate}
\item there exist a lower semicontinuous function $f : \D \to \R_+$ and a measurable function $g : \D\to\R_+$ such that
\begin{eqnarray*}
(i) & & u(z) \leq f(z) \, g(z) \quad \text{for every } z\in\D , \\[2mm]
(ii) & & \lim_{z\to e^{i\varphi} \atop  z \in \Omega_{\theta}(\varphi)} f(z) = 0  \text{ for every } \varphi\in [0,2\pi ) \text{ and some } \theta \in (0,\frac{\pi}{2}) , \\
(iii) & & \sup_{r\in (0,1)} \int_0^{2\pi } g(r e^{i\varphi} ) \; d\varphi <\infty ,
\end{eqnarray*}
\item there exists $m\in (0,\frac{\pi}{\theta} )$ such that
$$
M_r (u) = \sup_{\varphi\in [0,2\pi ]} u(re^{i\varphi } ) = O ( (1-r)^{-m} ) , \quad r\to 1- .
$$
\end{enumerate}
Then $u=0$.
\end{theorem}

\subsection{Analytic continuation across a linear boundary} \label{secanal}

In this section, we consider the square
$$
R := \{ z\in\C : -1 < {\rm Re}\, z <1 ,\, -1< {\rm Im}\, z < 1 \} ,
$$
and we study the question whether an analytic function
$u:R\setminus \R \to \C$ admits an analytic continuation to the
whole square $R$. For a thorough discussion of this type of problems see  \cite{Be72}. One of the analytic extension criteria is
provided by the following  classical {\it edge-of-the-wedge}
theorem.

\begin{proposition}[Edge-of-the-wedge] \label{eotw}
Let $u: R\setminus\R \to \C$ be an analytic function. Assume that
$$
\lim_{\beta\to 0} \big( u(\cdot + i\beta ) - u(\cdot -i\beta ) \big) = 0 
$$
on $(-1,1)$ in the sense of distributions. Then $u$ admits an
analytic extension to $R$.
\end{proposition}

\begin{remark}
Usually, the analytic extendability of $u$ is derived from the
existence and coincidence of the distributional limits
$\lim_{\beta\to 0+} u(\cdot + i\beta )$ and $ \lim_{\beta\to 0+}
u(\cdot -i\beta ).$ The (formally) more general Proposition
\ref{eotw} follows from \cite[Theorem C]{Ru71} using Carleman's trick as in \cite[Theorem A]{Wo47}; compare also with the proof of Theorem \ref{analytic} below.
\end{remark}

In particular, the conclusion of Proposition \ref{eotw} holds if
$\lim_{\beta\to 0} \big( u(\cdot + i\beta ) - u(\cdot -i\beta )
\big) = 0$ in $L^p$. However, if $\lim_{\beta\to 0} \big( u(\cdot
+ i\beta ) - u(\cdot -i\beta ) \big) = 0$ pointwise everywhere on
$(-1,1)$ then Proposition \ref{eotw} can hardly be applied
directly.

As a corollary to Theorem \ref{main2}, we obtain the following edge-of-the-wedge theorem
where distributional convergence is replaced by a combination of pointwise convergence everywhere and mean convergence. It improves the corresponding result in \cite[Theorem 3.1]{ChTo04}.

\begin{theorem} \label{analytic}
Let $u: R\setminus\R \to \C$ be an analytic function. Assume that
\begin{enumerate}
\item \label{2ana} there exist a lower semicontinuous function $f : R\setminus\R \to \R_+$ and a measurable function $g : R\setminus\R \to\R_+$ such that
\begin{eqnarray*}
(i) & & |u(z) - u(\bar{z}) | \leq f(z) \, g(z) \quad \text{for every } z\in R\setminus\R , \\[2mm]
(ii) & & \lim_{z\to \alpha \atop z\in \alpha +\Sigma_{\theta} } f(z) = 0 \mbox{ for every } \alpha\in (-1,1) \mbox{ and some } \theta \in (0,\frac{\pi}{2} ) , \\
(iii) & & \sup_{\beta \in (0,1)} \int_{-1}^1 g(\alpha + i\beta ) \; d\alpha <\infty .
\end{eqnarray*}
\item \label{1ana} there exists $m\geq 0$ such that
$$
\sup_{\alpha\in [-1,1 ]} | u(\alpha +i\beta ) - u(\alpha -i\beta ) | = O ( \beta^{-m} ) , \quad \beta\to 0+ ,
$$
\end{enumerate}

Then the function $u$ extends analytically to $R$.
\end{theorem}

\begin{proof}
The function $U(z):= (u(z ) - u(\bar{z}))/2$ is harmonic in $R\setminus \R.$ By the assumptions on $u$ and Theorem \ref{main2}, $U$ extends continuously to the whole rectangle $R$ and $U=0$ on the interval $(-1,1)$. We denote this extended function also by $U$. Since $U(z) = -U(\bar{z})$ and by the Schwartz reflection principle, $U$ is harmonic on $R$. 

There exists a harmonic conjugate $V$ of $U$ in $R,$ so that $\tilde{u} := U + i V$ is analytic in $R.$ The harmonic conjugate $V$ is unique up to an additive constant. 

The function $V_0 (z):=(u(z ) + u(\bar{z}))/(2i)$ is a harmonic conjugate of $U$ in $R\setminus\R$, and since $U(z)+iV_0(z) = u(z)$, 
it follows that $\tilde u - u$ is constant in each component of $R\setminus\R$. 
Since $V$ is continuous in $R$ and $V_0$ is symmetric with respect to $\R$, we obtain that
$\tilde{u}-u$ is constant in $R\setminus\R$. Thus, the function $u$ extends
analytically to $R.$
\end{proof}

\subsection{Stability of bounded $C_0$-semigroups}\label{semstab}

Let $X$ be a complex Banach space. Consider an abstract Cauchy
problem
\begin{equation} \label{acp}
\left\{ \begin{array}{ll}
\dot{u} (t) = A u(t) , & t\geq 0 ,\\[2mm]
u(0) = x , & x\in X,
\end{array} \right.
\end{equation}
where $A$ is a closed linear operator on $X$ with dense domain
$D(A).$  It is a fundamental fact of the theory of
$C_0$-semigroups that \eqref{acp} is well-posed if and only if $A$
generates a $C_0$-semigroup $\TT,$ \cite[Theorem 3.1.12]{ABHN01}. This semigroup is said to be
{\em stable} if for every $x \in X$ we have
$$
\lim_{t\to\infty} \| T(t)x \| = 0 .
$$

Recall that a function $u \in C(\mathbb R_+)$ is a {\it mild}
solution of \eqref{acp} if $\int_{0}^{t}u(s)\, ds \in D(A)$ and $A
\int_{0}^{t} u(s) \, ds = u(t) -x$ for all $t \ge 0.$ The
importance of the concept of stability comes from the fact that
the mild solutions of a well-posed abstract Cauchy problem
\eqref{acp} are precisely the orbits of the $C_0$-semigroup $\TT$
generated by $A.$ Thus, in the terminology of differential
equations, a stable $C_0$-semigroup corresponds to a well-posed
Cauchy problem \eqref{acp} for which all mild solutions are
asymptotically stable. One of the central problems in the theory
of stability is to characterize stability of a semigroup in {\em a
priori} terms of the generator, e.g. in terms of the spectrum of
the generator or, more generally, its resolvent.

For  accounts of (mostly spectral) stability theory, one may consult \cite{ABHN01}, \cite{EnNa99}, \cite{Ne96}. A discussion of recent developments can be found in \cite{ChTo06}. In this section, we present several resolvent stability criteria and their ``algebraic" counterparts. In a certain sense, these criteria are
optimal.

By the uniform boundedness principle, a stable semigroup is
automatically bounded, thus we will study stability of {\it bounded}
semigroups in the rest of the paper.

The following stability criterion from \cite[Theorem 3.1]{BaBrGr96} (see also \cite[Theorem 2.3]{Vu93}) will be fundamental for our study. Recall that a function $F:\R\to X$ is a {\em
complete trajectory} for a semigroup $\TT$ if $F(t+s) = T(t) F(s)$ for all $t\geq 0$, $s\in\R$. If $\TT$ is the adjoint of a $C_0$-semigroup, then $F$ is weak$^*$ continuous on $\mathbb R.$

\begin{theorem} \label{complete}
For a bounded $C_0$-semigroup $\TT$ on a Banach space $X$ the
following statements are equivalent:
\begin{itemize}
\item[(i)] The semigroup $\TT$ is stable,
\item[(ii)] The adjoint semigroup $(T(t)^*)_{t\geq 0}$ does not admit a bounded, nontrivial complete trajectory.
\end{itemize}
\end{theorem}

Theorem \ref{complete} can be considered as a topological characterization of stability. To relate this characterization to analytic properties of the generator we will need the notions of Carleman transform and Carleman spectrum.

For every bounded weakly measurable (weak$^*$ measurable, if $X$
is a dual space) function $F:\R \to X$  we define the {\em
Carleman transform} $\hat{F}$ by
$$
\hat{F} (\lambda ) := \left\{ \begin{array}{ll}
\int_0^\infty e^{-\lambda t} F(t) \; dt , & {\rm Re}\, \lambda >0 ,\\[2mm]
-\int_{-\infty}^0 e^{-\lambda t} F(t) \; dt , & {\rm Re}\, \lambda
<0 .
\end{array} \right.
$$
If $\TT$ is  a bounded $C_0$-semigroup on $X$ and $F(t)=T(t)x$ for
$t \ge 0,$ then
\begin{equation}\label{laplace}
 \hat{F} (\lambda )= R(\lambda, A), \qquad {\rm Re} \, \lambda > 0,
\end{equation}
where $R(\lambda, A)$ is the resolvent of the generator $A$ of
$\TT.$ \\
The Carleman transform $\hat{F}$ is analytic in
$\C\backslash i\R$. The set of singular points of $\hat{F}$ on
$i\mathbb R$, i.e. the set of all points near which $\hat{F}$ does
{\em not} admit an analytic extension, is called the {\it Carleman
spectrum} of $F.$

The Carleman spectrum of a bounded  function is nonempty unless the function is zero, see
for example \cite[Theorem 4.8.2]{ABHN01}; 
indeed, the only entire function $f$ satisfying $\| f(z)\|\le c/|{\rm Re}\, z|$, $z \in \mathbb C \setminus i\R$, is $0$.
Thus, by the equivalence
(i)$\Leftrightarrow$(ii) of Theorem \ref{complete}, a bounded $C_0$-semigroup is stable if and
only if for every bounded complete trajectory $F$ of the adjoint semigroup the Carleman transform
$\hat F$ extends to an entire function. This fact allows us to apply the analytic extension criterion
from Theorem \ref{analytic} to the study of stability of semigroups.

The following version of the resolvent identity for the Carleman transform of a complete trajectory appears to be useful for locating the Carleman spectrum of the trajectory, \cite[Lemma 6.1]{BaChTo02}.

\begin{lemma}\label{locresid}
Let $\TT$ be a bounded $C_0$-semigroup on a Banach space $X$ with
generator $A$, let $F$ be a bounded complete trajectory for
$(T(t)^*)_{t\geq 0},$ and let $\hat{F}$ be its Carleman transform.
Then for every $\lambda\in\C_+$ and
every $\mu\in\C\backslash i\R$
\begin{equation} \label{ri}
\hat{F}  (\mu ) = R(\lambda ,A^*)F(0) + (\lambda -\mu) R(\lambda
,A^*) \hat{F} (\mu ).
\end{equation}
\end{lemma}

In order to formulate our main results in this section, we need to recall that a Banach space $X$ has {\em Fourier type} $p\in [1,2]$, if the Fourier transform $\cF$ defined on the vector-valued Schwartz space $\cS (\R ;X)$ by
$$
\left(\cF \varphi \right) (\beta ) := \int_\R e^{-i\beta t}
\varphi (t) \; dt , \quad \beta\in\R,
$$
extends to a bounded linear operator from $L^p (\R ;X)$ into $L^q
(\R ;X),$  $\frac{1}{p}+\frac{1}{q}=1$, i.e. if the
Hausdorff-Young inequality holds for $X$-valued $L^p$-functions.
The class of Banach spaces with Fourier type greater than $1$ is
fairly large: it contains, for example, $L^p$-spaces with $p> 1$.
The Hilbert spaces can be characterized by the property of having
Fourier type $2.$ For a comprehensive survey of the theory of
Fourier type one may consult \cite{GKKT98}.

The following result essentially improves the stability criterion in \cite[Theorem 4.2]{ChTo04}.

\begin{theorem} \label{c0stable}
Let $A$ be the generator of a bounded $C_0$-semigroup $\TT$ on a
Banach space $X$ with Fourier type $p\in (1,2]$. Assume that the
set
$$
M:= \{ x\in X : \lim_{\alpha\to 0+} \alpha^{\frac{p-1}{p}} R(\alpha +i\beta ,A)x = 0 \text{ for every } \beta\in\R \}
$$
is dense in $X$. Then $\TT$ is stable.
\end{theorem}

\begin{remark}
A partial case of Theorem \ref{c0stable} for $p=2$ (i.e. when $X$ is a Hilbert space) was proved in \cite{To01} by means of an involved operator-theoretic construction. The reasoning there does not extend to the more general situation considered here.
\end{remark}

\begin{proof}
Let $F:\R\to X^*$ be a bounded complete trajectory for the adjoint
semigroup $(T(t)^*)_{t\geq 0}$, i.e. a weak$^*$ continuous function satisfying
\newline\nobreak
$T(t)^* F(s) = F(t+s)$ for every $t\geq 0$, $s\in\R$. Let $x\in M$
and let $F_x:=\langle F ,x\rangle$. Let $\hat{F_x}$  be the
Carleman transform of $F_x.$

By Lemma \ref{locresid}  for every $\alpha
>0$ and every $\beta\in\R$,
\begin{gather*}
|\hat{F_x} (\alpha +i\beta ) - \hat{F_x} (-\alpha +i\beta ) | =
| 2 \langle \alpha^{\frac{1}{p}} \hat{F_x} (-\alpha +i\beta ) , \alpha^{\frac{1}{q}} R(\alpha +i\beta ,A) x \rangle |\\
\le   g (\alpha +i\beta) \, f (\alpha +i\beta),
\end{gather*}
where
$$
g (\alpha +i\beta ):= \| 2 \alpha^{\frac{1}{p}} \hat{F_x} (-\alpha
+i\beta ) \|
$$
and
$$
f (\alpha +i\beta ):= \| \alpha^{\frac{1}{q}} R(\alpha +i\beta
,A)x \| .
$$

By the boundedness of $F$ and the Hausdorff-Young inequality,
\begin{equation} \label{Ggro}
\sup_{\alpha >0} \| g(\alpha +i\cdot )\|_{L^q (\R )} < \infty,
\quad {\textstyle \frac{1}{p} + \frac{1}{q} = 1 }.
\end{equation}
Moreover, by \eqref{laplace},
$$
\sup_{\alpha+i\beta \in\mathbb C_+}\alpha \|R(\alpha + i \beta, A)\| \le
\sup_{t \ge 0} \|T(t)\| < \infty.
$$
By assumption $x\in M$, and the resolvent identity implies for every
$\theta \in (0,\frac{\pi}{2})$ and every $\beta\in \R$ that
\begin{eqnarray*}
& & \limsup_{\alpha \to 0+ \atop |\beta' - \beta |\leq\alpha\,\tan\theta} f(\alpha +i\beta' ) \\
& \leq & \limsup_{\alpha \to 0+ \atop |\beta' - \beta |\leq\alpha\,\tan\theta} \| \alpha^{\frac{1}{q}} \bigl( R(\alpha +i\beta', A) - R(\alpha +i\beta , A) \bigr) x \| \\
& \leq & \limsup_{\alpha \to 0+ \atop |\beta' - \beta |\leq\alpha\,\tan\theta} \| \alpha\, \tan\theta \, R(\alpha +i\beta', A) \,  \alpha^{\frac{1}{q}} R(\alpha +i\beta, A) x \| \\
& \leq & \tan\theta \, \sup_{t\geq 0} \| T(t)\| \, \limsup_{\alpha \to 0+ } \| \alpha^{\frac{1}{q}} R(\alpha +i\beta , A) x \| \\
& = & 0 .
\end{eqnarray*}
It follows from this inequality, the boundedness of $F_x$ and
\eqref{Ggro} that we can apply Theorem \ref{analytic} in order to
see that the Carleman transform $\hat{F_x}$ extends analytically
through the imaginary axis to an entire function. This implies
$F_x=\langle F,x\rangle =0$. Since $M$ is dense in $X$, we conclude
that  $F=0$, i.e. there is no nontrivial bounded complete
trajectory for $(T(t)^*)_{t\geq 0}$. By Theorem \ref{complete}, this is
equivalent to the property that the semigroup $\TT$ is stable.
\end{proof}

\begin{remark}
Note that the conditions of Theorem \ref{c0stable} are not meaningful for $p=1.$ Indeed, $\lim_{\alpha\to 0+} R(\alpha +i\beta ,A)x = 0$ implies $x=0$ by the resolvent identity.
\end{remark}

If $A$ is the generator of a bounded $C_0$-semigroup, then there is a variety of (equivalent) ways to  define the {\em fractional powers} $(i\beta -A)^\gamma$ for every $\beta\in\R$ and every $\gamma >0$; see, for instance, \cite{MaSa01}. Note that $((i\beta
-A)^\gamma)_{\gamma
>0}$ form a semigroup and so the ranges ${\rm Rg}\, (i\beta
-A)^\gamma$ satisfy
$$
{\rm Rg}\, (i\beta -A)^\nu \subset {\rm Rg}\, (i\beta -A)^\gamma,
\quad  \beta\in\R,
$$
whenever $\nu > \gamma.$ We proceed with a corollary of Theorem
\ref{c0stable} which is a kind of algebraic criterion for
stability. It can be considered as a limit case of Theorem
\ref{c0stable}.

\begin{corollary} \label{c0stablecor}
Let $A$ be the generator of a bounded $C_0$-semigroup 
\newline\nobreak$\TT$ on a
Banach space $X$ with Fourier type $p\in [1,2]$. If
\begin{equation} \label{range}
\bigcap_{\beta \in \mathbb R} {\rm Rg}\, (i\beta -A)^{\frac{1}{p}}
\mbox{ is dense in } X ,
\end{equation}
then $\TT$ is stable.
\end{corollary}

\begin{remark}
Corollary \ref{c0stablecor} improves the stability result from \cite[Corollary 4.6]{ChTo04}; the exponent $\frac{1}{p}$ in \eqref{range} was there replaced by an exponent $\gamma >\frac{1}{p}$.

Corollary \ref{c0stablecor} is optimal in the following sense: for every $p\in [1,2]$ there exist a Banach space $X$ having Fourier type $p$ and a {\em nonstable} bounded semigroup with generator $A$ such that for {\em every} $\gamma\in (0,\frac{1}{p})$
$$
\bigcap_{\beta \in \mathbb R} {\rm Rg}\, (i\beta -A)^{\gamma}
\mbox{ is dense in } X ;
$$
see \cite[Example 4.9]{ChTo03}.

Note that the stability condition from Corollary \ref{c0stablecor} is sufficient but in general not necessary for stability, except possibly in the case $p=2$, \cite[Example 4.1]{ChTo03}. In fact, it is an open problem whether the condition \eqref{range} (or the stability condition in Theorem \ref{c0stable}) characterizes stability in Banach spaces with Fourier type $p=2$, i.e in Hilbert spaces.
\end{remark}

\begin{proof}[Proof of Corollary \ref{c0stablecor}]
The case $p=1$ has been proved in \cite[Theorem 2.4]{BaChTo02}. The assumption implies that ${\rm Rg}\, (i\beta -A)$ is dense in $X$. By the fractional mean ergodic theorem, \cite[Proposition 2.3]{Ko69IV}, \cite[Proposition 2.2]{Ws98}, for every $x\in X$ and every $p>1$
$$
\lim_{\alpha\to 0+} \alpha^{\frac{p-1}{p}} R(\alpha +i\beta ,A)^{\frac{p-1}{p}} x = 0 .
$$
By \cite[Theorem 2.4]{Ko69IV}, \cite[Theorem 2.4]{Ws98}, for every $x\in {\rm Rg}\, (i\beta -A)^{\frac{1}{p}}$
$$
\lim_{\alpha\to 0+} R(\alpha +i\beta ,A)^{\frac{1}{p}} x \text{ exists}.
$$
Hence, if $x\in {\rm Rg}\, (i\beta -A)^{\frac{1}{p}}$ and $p>1$, then
$$
\lim_{\alpha\to 0+} \alpha^{\frac{p-1}{p}} R(\alpha +i\beta ,A) x = \lim_{\alpha\to 0+} \alpha^{\frac{p-1}{p}} R(\alpha +i\beta ,A)^{\frac{p-1}{p}} R(\alpha +i\beta ,A)^{\frac{1}{p}} x = 0 .
$$
The claim for $p>1$ thus follows from Theorem \ref{c0stable}.
\end{proof}

\begin{remark}
We conclude with a remark concerning stability of discrete operator semigroups $(T^n)_{n \ge 0}$, where $T$ is a bounded linear operator on $X.$ {\em Stability} for a bounded linear operator $T$ means that $\lim_{n\to\infty} \| T^n x\| =0$ for every $x\in X$. All the above arguments can be directly transferred to the discrete case after an appropriate change of notions. Sketches of this procedure can be found in \cite{ChTo04} and in \cite{BaChTo02}. We obtain the following analog 
of Theorem \ref{c0stable} in the discrete setting.

\begin{theorem}
Let $T$ be a power bounded linear operator  on a Banach space $X$
with Fourier type $p\in (1,2]$. Assume that the set
$$
\{ x \in X: \lim_{r \to 1+}  \| (r-1)^{\frac{p-1}{p}} R(r \xi, T)x
\| = 0 \text{ for every } \xi \in \mathbb T \}
$$
is dense in $X.$ Then $T$ is stable.
\end{theorem}

Also Corollary \ref{c0stablecor} admits a discrete counterpart.

\begin{theorem}
Let $T$ be a power bounded linear operator on a Banach space $X$ with Fourier type $p\in [1,2]$. If
$$
\bigcap_{\xi \in \mathbb T} {\rm Rg}\, (\xi -T)^{\frac{1}{p}} \mbox{ is dense in } X ,
$$
then $T$ is stable.
\end{theorem}

\end{remark}

\section{Appendix} \label{appe}

In this Appendix, we show that Wolf's theorem (Theorem \ref{wolf}) is optimal with respect to the growth condition on $M_r (u)$. 

\begin{example} Given $\theta\in (0,\pi/2)$ and $\varepsilon_0>0$, there exists 
a function $u\ne 0$ harmonic in the unit disc such that
\begin{eqnarray*}
(i) & & \lim_{z\to e^{i\varphi}
\atop z \in \Omega_{\theta}(\varphi)} u(z) = 0 \mbox{ for every } \varphi \in [0,2\pi );\\
(ii) & &
M_r(u)={\rm O} ( e^{\varepsilon_0 (1-r)^{-\frac{\pi}{2\theta}}} ) \quad
\text{ as } r\to 1-.
\end{eqnarray*}\label{app}
\end{example}

For a similar construction see \cite[Example~7.1]{BoVo89}.

\begin{proof}[Sketch of the proof]
Choose small $\beta>0$ and $\varepsilon>0$. Let $D_0$ be the union of $\mathbb D$ and a small neighborhood (in $\mathbb C$) of the arc
$$
\bigl\{z \in \partial\mathbb D\setminus\{1\}: \frac\pi2-\theta-\frac{2\beta}3\le\arg(1-z)\le\frac\pi2-\theta-\frac{\beta}3\bigr\},
$$
and let 
$$
\Delta=\bigl\{z \in \overline{D_0}\setminus\{1\}: \arg(1-z)>\frac\pi2-\theta-\beta\bigr\}.
$$
We choose $K\in\mathbb R$ such that if $\arg(1-z)=\frac \pi2-\theta$, then
$K+i\varepsilon\frac{1+z}{1-z}\in e^{i\theta}\mathbb R_+$.
For large even $A$ we define
$$
f(z)= \left\{ \begin{array}{ll}
\bigl( \frac{1-z}{1+z} \bigr)^A \exp \Bigl( \bigl( K+i\varepsilon\frac{1+z}{1-z} \bigr)^{\pi/(2\theta)} \Bigr) & \text{if } z\in \Delta, \\[2mm]
0  & \text{if } z\in \overline{D_0}\setminus\Delta .
              \end{array}
\right.
$$
Then $\Im f(z)=0$ for every $z\in\Delta\cap\partial\mathbb D$, 
and $|f(z)|+|\nabla f(z)|=O(|1-z|^2)$ as $z\to 1$, $\frac\pi2-\theta-\beta<\arg(1-z)\le\frac \pi2-\theta$.

Next, we fix $\tilde g\in C^2 ([-\pi/2,\pi/2])$ such that $\tilde g{\bigm|} [\frac\pi2-\theta-\frac{\beta}{3},\frac\pi2]=1$, 
$\tilde g{\bigm|} [-\frac\pi2,\frac\pi2-\theta-\frac{2\beta}{3}]=0$, $0\leq \tilde{g}\leq 1$, and we define $g(z):=\tilde g(\arg(1-z))$, $z\in \overline{D_0}\setminus\{1\}$.
Then $g\in C^2(\overline{D_0}\setminus\{1\})$ and 
$$
|\nabla g(z)|\le \frac{C}{|1-z|}, \qquad z\in D_0.
$$

Set $h:=fg$. Then $h\in C^2(D_0)\cap C(\overline{D_0}\setminus\{1\})$, 
$\Im h(z)=0$ for every $z\in \partial\mathbb D\setminus\{1\}$,
\begin{eqnarray*}
& & |\bar\partial h(z)|={\rm O}\,(|1-z|) \text{ as } z\to 1, \text{and } \\
& & \lim_{z\to 1 \atop z\in\Omega_\theta(0)} h(z) =  0 . 
\end{eqnarray*}
Next, we set
$$
h_1(z)=\frac{1}{2\pi}\int_{D_0}\frac{\bar\partial h(\zeta)}{z-\zeta}dm_2(\zeta).
$$
Then $h_1\in C^1(D_0)\cap C(\overline{\mathbb D})$ and $\bar\partial h_1=\bar\partial h$ on $D_0$. Therefore, the function $H=h-h_1$ is analytic in $D_0$, $H\in C(\overline{\mathbb D}\setminus\{1\})$, and
$$
\log|H(z)|\le C+\frac{C\varepsilon}{|1-z|^{\pi/(2\theta)}},\qquad z\in \mathbb D.
$$
Let $v$ be the Poisson extension of ${\rm Im}\, h_1{\bigm|}\partial\mathbb D$ to $\mathbb D$, and let $w:=v+\Im H$. Then $w$ is harmonic in
$\mathbb D$, $w\in C(\overline{\mathbb D}\setminus\{1\})$, $w=0$ on $\partial\mathbb D\setminus\{1\}$, and
\begin{eqnarray*}
\log|w(z)| & \le & C+\frac{C\varepsilon}{|1-z|^{\pi/(2\theta)}} \\
& \leq &  C+\frac{C\varepsilon}{(1-|z|)^{\pi/(2\theta)}} ,\qquad z\in \mathbb D .
\end{eqnarray*}
Finally,
$$
\lim_{\stackrel{z\to 1}{z\in\Omega_\theta(0)}}w(z)= \lim_{\stackrel{z\to 1}{z\in\Omega_\theta(0)}}(v(z) - {\Im h_1}(z))  + \lim_{\stackrel{z\to 1}{z\in\Omega_\theta(0)}} {\Im h} (z) =0.
$$
\end{proof}

 \providecommand{\bysame}{\leavevmode\hbox to3em{\hrulefill}\thinspace}

\bibliographystyle{amsplain}

\providecommand{\bysame}{\leavevmode\hbox to3em{\hrulefill}\thinspace}
\providecommand{\MR}{\relax\ifhmode\unskip\space\fi MR }
\providecommand{\MRhref}[2]{%
  \href{http://www.ams.org/mathscinet-getitem?mr=#1}{#2}
}
\providecommand{\href}[2]{#2}

\end{document}